\documentclass[12pt]{amsart}
\usepackage[cp1251]{inputenc}
\usepackage[english]{babel}
\usepackage{amsmath}
\usepackage{amssymb}
\usepackage{amsfonts}

\begin{document}


\title[On discrete boundary value problems]
      {On discrete boundary value problems with nonlocal conditions in a quarter-plane}
\author{Vladimir Vasilyev}
\address{Chair of Applied Mathematics and Computer Modeling\\
    Belgorod State National Research University\\
         Pobedy street 85, Belgorod 308015, Russia}

        \email{vbv57@inbox.ru}

  \author{Anastasia Mashinets}
\address{Chair of Applied Mathematics and Computer Modeling\\
    Belgorod State National Research University\\
         Pobedy street 85, Belgorod 308015, Russia}

        \email{anastasia.kho@yandex.ru}


\begin{abstract}
We consider discrete analogue of model pseudo-differential equations in discrete plane sector using discrete variant of Sobolev--Slobodetskii spaces. Starting from the concept of wave factorization for elliptic periodic symbol we describe solvability conditions for the equations and corresponding discrete boundary value problems. We give also a comparison between discrete and continuous solutions in appropriate discrete normed spaces.
\end{abstract}

\maketitle

\section{Introduction}

A theory of pseudo-differential operators and equations \cite{H83,Ta81,Tr80,E81} not so long history than other mathematical subjects of analysis. Nevertheless, these operators and related boundary value problems widely arise in a lot of applied problems in physics and technique (see, for example \cite{SFO} and references therein). Discrete aspects of the theory are reflected in mathematical papers more weak \cite{R,Ru} although these studies are closely related to the theory of Fourier series \cite{Ed82}. In our opinion the discrete theory is very important since it permits to use computer calculations to solve concrete applied problems.

We interested in studying discrete pseudo-differential equations and their solvability in appropriate discrete functional spaces. There are certain approaches to studying discrete boundary value problems for partial differential equations including finite difference method \cite{S01,R02}. But these approaches are not applicable to studying discrete boundary value problems for elliptic pseudo-differential equations. According to this statement the first author with colleagues  has started to develop discrete theory for elliptic pseudo-differential equations \cite{VV6,TV}. This is main motivation, and we have started from certain canonical domains. First considerations were related to discrete $m$-dimensional space and half-space, and here we consider discrete quadrant.

We consider a special type of boundary conditions, namely integral conditions on a boundary. These conditions are nonlocal, and it seems, such conditions are artificial. But there are a lot of applied problems for partial differential equations with such boundary conditions \cite{AB,CK,KMT}, therefore it is natural way. Moreover, these conditions appear in a natural way to determine arbitrary functions in a general solution of an elliptic pseudo-differential equation.

\section{Cones, periodic symbols, digital operators and equations}

\subsection{Discrete spaces and transforms}
Let $\mathbb Z^2$ be an integer lattice in a plane. Let $K=\{x\in\mathbb R^2: x=(x_1,x_2), x_1>0, x_2>0\}$ be a quadrant, $K_d=h\mathbb Z^2\cap K, h>0$. We consider functions of discrete variable $u_d(\tilde x), \tilde x=(\tilde x_1,\tilde x_2)\in h\mathbb Z^2$.

Let us denote $\mathbb T^2=[-\pi,\pi]^2, \hbar=h^{-1}$. We consider functions defined in $\hbar\mathbb T^2$ as periodic functions defined in $\mathbb R^2$ with basic square of periods $\hbar\mathbb T^2$.

One can define the discrete Fourier transform for the function $u_d$
\[
(F_du_d)(\xi)\equiv\tilde u_d(\xi)=\sum\limits_{\tilde x\in h\mathbb Z^2}e^{-i\tilde x\cdot\xi}u_d(\tilde x)h^2,~~~\xi\in\hbar\mathbb T^2,
\]
if the latter series converges, and the function   $\tilde u_d(\xi)$ is a periodic function in $\mathbb R^2$ with basic square of periods $\hbar\mathbb T^2$. Such discrete Fourier transform preserves all properties of integral Fourier transform, and the inverse discrete Fourier transform looks as follows
 \[
 (F_d^{-1}\tilde u_d)(\tilde x)=\frac{1}{(2\pi)^2}\int\limits_{\hbar\mathbb T^2}e^{i\tilde x\cdot\xi}\tilde u_d(\xi)d\xi,~~~\tilde x\in h\mathbb Z^2.
 \]

The discrete Fourier transform gives one-to-one correspondence between spaces  $L_2(h\mathbb Z^2)$ and $L_2(\hbar\mathbb T^2)$ with norms
\[
||u_d||_2=\left(\sum\limits_{\tilde x\in h\mathbb Z^2}|u_d(\tilde x)|^2h^2\right)^{1/2},
\quad
||\tilde u_d||_2=\left(\int\limits_{\xi\in\hbar\mathbb T^2}|\tilde u_d(\xi)|^2d\xi\right)^{1/2}.
\]

We need more general discrete functional spaces and we introduce such spaces using divided differences \cite{S01}.

The divided differences of first order look as follows
\[
(\Delta_1^{(1)}u_d)(\tilde x)=h^{-1}(u_d(\tilde x_1+h,\tilde x_2)-u_d(\tilde x_1,\tilde x_2)),
\]
\[
(\Delta_2^{(1)}u_d)(\tilde x)=h^{-1}(u_d(\tilde x_1,\tilde x_2+h)-u_d(\tilde x_1,\tilde x_2)),
\]
and their discrete Fourier transforms are given by formulas
\[
\widetilde{(\Delta_k^{(1)}u_d)}(\xi)=h^{-1}(e^{-ih\cdot\xi_k}-1)\tilde u_d(\xi), k=1,2.
\]

The divided difference of second order is a divided difference of first order from divided difference of first order
\[
(\Delta_1^{(2)}u_d)(\tilde x)=h^{-2}(u_d(\tilde x_1+2h,\tilde x_2)
-2u_d(\tilde x_1+h,\tilde x_2)+u_d(\tilde x_1,\tilde x_2)),
\]
\[
(\Delta_2^{(2)}u_d)(\tilde x)=h^{-2}(u_d(\tilde x_1,\tilde x_2+2h)
-2u_d(\tilde x_1,\tilde x_2+h)+u_d(\tilde x_1,\tilde x_2)),
\]
with the Fourier transform
\[
\widetilde{(\Delta_k^{(2)}u_d)}(\xi)=h^{-2}(e^{-ih\cdot\xi_k}-1)^2\tilde u_d(\xi), k=1,2.
\]

Discrete analogue of the Laplacian is the following
\[
(\Delta_du_d)(\tilde x)=(\Delta_1^{(2)}u_d)(\tilde x)+(\Delta_2^{(2)}u_d)(\tilde x),
\]
so that its Fourier transform is
\[
\widetilde{(\Delta_du_d)}(\xi)=h^{-2}((e^{-ih\cdot\xi_1}-1)^2+(e^{-ih\cdot\xi_2}-1)^2)\tilde u_d(\xi).
\]

We use such discrete objects for constructing discrete Sobolev--Slobodetskii spaces to study wide class of discrete equations.

First, we introduce discrete analogue of the Schwartz space
 $S(h\mathbb Z^2)$ as a set of discrete functions with finite semi-norms
 \[
 |u_d|=\sup\limits_{\tilde x\in h\mathbb Z^2}(1+|\tilde x|)^l|\Delta^{({\bf k})}u_d(\tilde x)|
 \]
 for arbitrary $l\in\mathbb N, {\bf k}=(k_1,k_2), k_r\in\mathbb N, r=1,2$,
 \[
 \Delta^{({\bf k})}u_d(\tilde x)=\Delta^{k_1}_1\Delta^{k_2}_2u_d(\tilde x).
 \]

 {\bf Definition 1.} {\it
 A discrete distribution is called a linear continuous functional defined on the space  $S(h\mathbb Z^2)$.
 }

 A set of such distributions will be denoted by  $S'(h\mathbb Z^2)$, and a value of the discrete distribution $f_d$ on the test discrete function  $u_d\in S(h\mathbb Z^2)$ will be denoted by $(f_d,u_d)$.

 One can introduce a concept of a support for a discrete distribution. Namely, a support of the discrete function $u_d\in S(h\mathbb Z^2)$ is a subset of the set  $h\mathbb Z^2$ such that  $u_d(\tilde x)\neq 0$ for all points $\tilde x$ from this subset. For an arbitrary set  $M\subset\mathbb R^2$ we denote $M_d=M\cap h\mathbb Z^2$, and then one says that  $f_d=0$
in the discrete domain $M_d$ if $(f_d,u_d)=0, \forall u_d\in S(M_d),$ where $S(M_d)\subset S(h\mathbb Z^2)$ consists of discrete functions with supports in  $M_d$. If $\widetilde M_d$ is a union of such $M_d$ where $f_d=0$ then support of the discrete distribution  $f_d$ is the set $h\mathbb Z^2\setminus\widetilde M_d$.

Similarly \cite{Vl} we can define standard operations in the space $S'(h\mathbb Z^2)$, but differentiation will be changed by divided difference of first order. These operations are described in  \cite{VV6} in details, a convergence is meant as a weak convergence in the space  $S'(h\mathbb Z^2)$.

{\bf Example 1.}
 If the function $f_d(\tilde x)$ is locally summable then it generates the discrete distribution
 \begin{equation}\label{1}
 (f_d,u_d)=\sum\limits_{\tilde x\in h\mathbb Z^2}f_d(\tilde x)u_d(\tilde x)h^2,~~~\forall u_d\in S(h\mathbb Z^2).
 \end{equation}

But there are different possibilities, for example, analogue of the Dirac mass-function
\[
(\delta_d,u_d)=u_d(0),
\]
which can not be represented by the formula  \eqref{1}.

Let $\zeta^2=h^{-2}((e^{-ih\cdot\xi_1}-1)^2+(e^{-ih\cdot\xi_2}-1)^2)$. We introduce the following definition.

{\bf Definition 2.} {\it
The space $H^s(h\mathbb Z^2)$ consists of discrete distributions and it is a closure of the space $S(h\mathbb Z^2)$ with respect to the norm
\begin{equation}\label{2}
||u_d||_s=\left(\int\limits_{\hbar\mathbb T^2}(1+|\zeta^2|)^s|\tilde u_d(\xi)|^2d\xi\right)^{1/2}.
\end{equation}
}

Let us remind that a lot of properties of such discrete spaces were studied in  \cite{F}, Varying the parameter $h$ in \eqref{2} we obtain different norms which are equivalent to the  $L_2$-norm. But constants in this equivalence depend on  $h$. In our constructions all constants do not depend on $h$.

{\bf Definition 3.} {\t
The space $H^s(K_d)$ consists of discrete distributions from  $H^s(h\mathbb Z^2)$ such that their supports belong to the set  $\overline{K_d}$. A norm in the space  $H^s(K_d)$ is induced by the norm of the space   $H^s(h\mathbb Z^2)$. The space $H^s_0(K_d)$ consists of discrete distributions   $f_d\in S'(h\mathbb R^2)$ with supports inside of $K_d$, and these discrete distributions must admit a continuation into the space  $H^s(h\mathbb Z^2)$. A norm in the space  $H^s_0(K_d)$ is given by the formula
\[
||f_d||^+_s=\inf||\ell f_d||_s,
\]
where infimum is taken for all continuations $\ell$.
}

The Fourier image of the space $H^s(K_d)$ will be denoted by $\widetilde H^s(K_d)$.

\subsection{Symbols, operators and projectors}
Let $\widetilde A_d(\xi)$ be a measurable periodic function in $\mathbb R^2$ with basic square of periods $\hbar\mathbb T^2$. Such functions we call symbols.

{\bf Definition 4.} {\it
A digital pseudo-differential operator  $A_d$ with the symbol $A_d(\xi)$ in the discrete quadrant $K_d$ is called an operator of the following type
\begin{equation}\label{3}
(A_du_d)(\tilde x)=\sum\limits_{\tilde y\in h\mathbb Z^2}h^2\int\limits_{\hbar\mathbb T^2}\widetilde A_d(\xi)e^{i(\tilde x-\tilde y)\cdot\xi}\tilde u_d(\xi)d\xi,~~~\tilde x\in K_d,
\end{equation}
}

We say that the operator $A_d$ is elliptic one if
\[
ess~\inf_{\xi\in\hbar\mathbb T^2}|A_d(\xi)|>0.
\]

A more general digital pseudo-differential operator with the symbol   $\widetilde A_d(\tilde x,\xi)$ depending on a spatial variable  $\tilde x$
\[
(A_du_d)(\tilde x)=\sum\limits_{\tilde y\in h\mathbb Z^2}h^2\int\limits_{\hbar\mathbb T^2} A_d(\tilde x,\xi)e^{i(\tilde x-\tilde y)\cdot\xi}\tilde u_d(\xi)d\xi,~~~\tilde x\in K_d,
\]
can be defined in the same way, but here we consider only operators of type \eqref{3}.

We consider symbols satisfying the condition
\begin{equation}\label{100}
c_1(1+|\zeta^2|)^{\alpha/2}\leq|A_d(\xi)|\leq c_2(1+|\zeta^2|)^{\alpha/2}
\end{equation}
with constants $c_1, c_2$ non-depending on $h$. The number
 $\alpha\in\mathbb R$ is called an order of digital pseudo-differential operator $A_d$.

 The following simple result can be proved easily.

{\bf Lemma 1.} {\it
 A digital pseudo-differential operator  $A_d$ with the symbol $\widetilde A_d(\xi)$ is a linear bounded operator
$H^s(h\mathbb Z^2)\to H^{s-\alpha}(h\mathbb Z^2)$ with a norm non-depending on  $h$.
}

We study a solvability of the discrete equation
\begin{equation}\label{4}
(A_du_d)(\tilde x)=v_d(\tilde x),~~~\tilde x\in K_d,
 \end{equation}
 in the space $H^s(K_d)$ assuming that $v_d\in H^{s-\alpha}_0(K_d)$.

We will use certain special domain in two-dimensional complex space  $\mathbb C^2.$ A domain of the type ${\mathcal T}_h(K)=\hbar\mathbb T^2+iK$ is called a tube domain over the quadrant $K$, and we will consider analytical functions $f(x+i\tau)$ in the domain ${\mathcal T}_h(K)=\hbar\mathbb T^2+iK$.

Let us introduce the periodic Bochner kernel similar \cite{Vl}
\[
B_h(z)=\sum\limits_{\tilde x\in K_d}e^{i\tilde x\cdot(\xi+i\tau)}h^2,~~~\xi\in\hbar\mathbb T^2,~~~\tau\in K,
\]
and corresponding integral operator
\[
(B_h\tilde u_d)(\xi)=\lim\limits_{\tau\to 0, \tau\in K}\frac{1}{4\pi^2}\int\limits_{\hbar\mathbb T^2}B_h(\xi+i\tau-\eta)\tilde u_d(\eta)d\eta.
\]

{\bf Lemma 2.} {\it
For the quadrant $K$ the operator $B_h$ has the following form
\[
(B_h\tilde u_d)(\xi)=\frac{h^2}{8\pi^2}\int\limits_{\mathbb T^2}\tilde u_d(\eta)d\eta+\lim\limits_{\tau\to 0+}\frac{ih}{8\pi^2}\int\limits_{\mathbb T^2}\cot\frac{h(\xi_1-\eta_1+i\tau_1)}{2}\tilde u_d(\eta)d\eta+
\]
\[
+\lim\limits_{\tau\to 0+}\frac{ih}{8\pi^2}\int\limits_{\mathbb T^2}\cot\frac{h(\xi_2-\eta_2+i\tau_2)}{2}\tilde u_d(\eta)d\eta-
\]
\[
-\lim\limits_{\tau\to 0+}\frac{h^2}{8\pi^2}\int\limits_{\mathbb T^2}\cot\frac{h(\xi_1-\eta_1+i\tau_1)}{2}\cot\frac{h(\xi_2-\eta_2+i\tau_2)}{2}\tilde u_d(\eta)d\eta,
\]
and $B_h$ is a linear bounded operator $H^s(\hbar\mathbb T^2)\rightarrow H^s(\hbar\mathbb T^2)$ for $|s|<1/2$. Moreover, the operator $B_h$ is a projector $\widetilde H^s(h\mathbb Z^2)\rightarrow \widetilde H^s(K_d)$.
}

{\bf Proof.}
Corresponding calculations for one-dimensional discrete cone were done in \cite{VV1}. We use these evaluations adapting to our two-dimensional case. Since
\[
\sum\limits_{\tilde x_k\in h\mathbb Z_+}e^{-i\tilde x_kz_k}h=\frac{h}{2}-\frac{ih}{2}\cot\frac{hz_k}{2},~~~z_k=\xi_k+i\tau_k,~~~k=1,2.
\]
then multiplying two factors and applying the Fourier property on correspondence between a product and convolution we obtain the assertion.

Boundedness of the one-dimensional operator with the kernel $h\cot\frac{hz}{2}$for $|s|<1/2$ was proved in \cite{F}, Theorem 6; two-dimensional case can be considered by the same method.
$\blacksquare$

{\bf Remark 1.}
The operator $B_h$ is so called periodic bi-singular operator. Using classical results for Cauchy type integral \cite{Ga,Mu} one can evaluate the boundary value, but it is not important this time.
Since these formulas are very huge we can do some simplifications without lost of generality. For example, we can consider the space
$S_1(h\mathbb Z^2)\subset S_1(h\mathbb Z^2)$ with zeroes in coordinate axes and consider the space $H^s(h\mathbb Z^2)$ as closure of the set
$S_1(h\mathbb Z^2)$ assuming that all functions of discrete variable vanish on coordinate axes.
For this case the first three summands in $B_h$ will be zero.

{\bf Lemma 3.} {\it If $|s|<1/2$ then the space $\widetilde H^s(h\mathbb Z^2)$ is uniquely represented as the direct sum
\[
\widetilde H^s(h\mathbb Z^2)=\widetilde H^{s}(K_d)\oplus\widetilde H^{s}(h\mathbb Z^2\setminus K_d)
\]
}

{\bf Proof.}
It is simple consequence of Lemma 2. Indeed, the unique representation of the function $\tilde f\in\widetilde H(h\mathbb Z^2)$ is the following
\[
\tilde f=B_h\tilde f+(I-B_h)\tilde f.
\]
A uniqueness of the such representation is possible only for $|s|<1/2$.
$\blacksquare$

To describe a solvability picture for the discrete equation \eqref{4} we need some additional elements of multidimensional complex analysis. We give it in the next section.

\section{Periodic wave factorization}

This concept is a periodic analogue of the wave factorization \cite{V0}. Some first preliminary considerations and results were described in \cite{V2,V3,V4,V5}.

{\bf Definition 5.} {\it
A periodic wave factorization for the elliptic symbol
$A_d(\xi)\in E_{\alpha}$ is called its representation in the form
\[
A_d(\xi)=A_{d,\neq}(\xi)A_{d,=}(\xi),
\]
where the factors  $A_{d,\neq}(\xi), A_{d,=}(\xi)$ admit analytical continuation into tube domains ${\mathcal T}_h(K), {\mathcal T}_h(-K)$ respectively with estimates
\[
c_1(1+|\hat\zeta^2|)^{\frac{\ae}{2}}\leq|A_{d,\neq}(\xi+i\tau)|\leq c'_1(1+|\hat\zeta^2|)^{\frac{\ae}{2}},
\]
\[
c_2(1+|\hat\zeta^2|)^{\frac{{\alpha-\ae}}{2}}\leq|A_{d,=}(\xi-i\tau)|\leq c'_2(1+|\hat\zeta^2|)^{\frac{{\alpha-\ae}}{2}},
\]
and constants $c_1, c'_1, c_2, c'_2$ non-depending on $h$, where
\[
\hat\zeta^2\equiv\hbar^2\left((e^{-ih(\xi_1+i\tau_1)}-1)^2+(e^{-ih(\xi_2+i\tau_2)}-1)^2\right),
\]
\[
\xi=(\xi_1,\xi_2)\in\hbar\mathbb T^2,~~~\tau-(\tau_1,\tau_2)\in K.
\]

The number $\ae\in{\mathbb R}$ is called an index of periodic wave factorization.
}

Unfortunately, we have no an algorithm to construct the factorization. But the are certain examples of periodic symbols which admit such factorization. We give one of them.

If $f$ is an arbitrary function of a discrete variable,  $f\in S(h\mathbb Z^2)$, $supp~f\subset K_d\cup(-K_d)$ then we have
\[
f=\chi_+f+\chi_-f,
\]
where $\chi_{\pm}$ are indicators of $\pm K_d$. Applying the discrete Fourier transform
we obtain the representation  $\tilde f=\tilde f_++\tilde f_-$, and $\tilde f_{\pm}$ admit an analytical continuation into  ${\mathcal T}_{h}(\pm K)$ according to Lemma 2. Thus, we can write $\exp\tilde f=\exp\tilde f_+\cdot\exp\tilde f_-$, therefore we obtain periodic wave factorization with index zero for the function $\exp\tilde f$.

{\it Everywhere below we assume existence of such periodic wave factorization for the symbol  $A_d(\xi)$  with index $\ae$.}

\subsection{A unique solvability}

This section is devoted to most simple case when a solution of the equation \eqref{4} exists and it is unique.

{\bf Theorem 1.} {\it
Let $|\ae-s|<1/2$. Then the equation  \eqref{4} has a unique solution for arbitrary right hand side $v_d\in H^{s-\alpha}_0(K_d)$, and it is given by the formula
\[
\tilde u_d(\xi)-A^{-1}_{d,\neq}(\xi)B_h(A^{-1}_{d,=}(\xi)\widetilde{(\ell v_d)}(\xi)),
\]
where $\ell v_d$ is an arbitrary continuation of $v_d$ into $H^{s-\alpha}(h\mathbb Z^2)$.
}

{\bf Proof}
Let $\ell v_d$ be an arbitrary continuation of $v_d\in H^{s-\alpha}_0(K_d)$ into $H^{s-\alpha}(h\mathbb Z^2)$. Let us introduce the function
\[
w_d(\tilde x)=(\ell v_d)(\tilde x)-(A_du_d)(\tilde x),
\]
so that $w(\tilde x)=0$ for $\tilde x\notin K_d$.

Now we write the equation \eqref{4} in the form
\[
(A_du_d)(\tilde x)+w_d(\tilde x)=(\ell v_d)(\tilde x),~~~\tilde x\in h\mathbb Z^2,
\]
and after applying the discrete Fourier transform and periodic wave factorization we obtain
\begin{equation}\label{5}
A_{d,\neq}(\xi)\tilde u_d(\xi)+A^{-1}_{d,=}(\xi)\tilde w_d(\xi)=A^{-1}_{d,=}(\xi)\widetilde{(\ell v_d)}(\xi),~~~\xi\in\hbar\mathbb T^2,
\end{equation}

We have the following inclusions according to Lemma 1 and Lemma 2
\[
A_{d,\neq}(\xi)\tilde u_d(\xi)\in\widetilde H^{s-\ae}(K_d),~~~A^{-1}_{d,=}(\xi)\tilde w_d(\xi)\in\widetilde H^{s-\ae}(h\mathbb Z^2\setminus K_d),
\]
\[
A^{-1}_{d,=}(\xi)\widetilde{(\ell v_d)}(\xi)\in\widetilde H^{s-\ae}(h\mathbb Z^2),
\]
and then according to Lemma 3 the right hand side of the equality \eqref{5} is uniquely represented by the sum
\[
A^{-1}_{d,=}(\xi)\widetilde{(\ell v_d)}(\xi)=f^+_d(\xi)+f^-_d(\xi),
\]
where
\[
f^+_d(\xi)=B_h(A^{-1}_{d,=}(\xi)\widetilde{(\ell v_d)}(\xi)),~~~f^-_d(\xi)=(I-B_h)(A^{-1}_{d,=}(\xi)\widetilde{(\ell v_d)}(\xi))
\]  .

Further, we rewrite the equality  \eqref{5}
\[
A_{d,\neq}(\xi)\tilde u_d(\xi)-f^+_d(\xi)=f^-_d(\xi)-A^{-1}_{d,=}(\xi)\tilde w_d(\xi)
\]
and using the uniqueness of the representation as the direct sum  $\widetilde H^{s-\ae}(K_d)\oplus\widetilde H^{s-\ae}(h\mathbb Z^2\setminus K_d)$ we conclude that both left hand side and right hand side should be zero. Thus,
\[
\tilde u_d(\xi)=A^{-1}_{d,\neq}(\xi)B_h(A^{-1}_{d,=}(\xi)\widetilde{(\ell v_d)}(\xi)),
\]
and Theorem 1 is proved.
$\blacksquare$

\section{Discrete boundary value problem}

In this section we consider more interesting case when the equation  \eqref{4} has a lot of solutions.

\subsection{Form of a discrete solution}

This section uses some results from  \cite{VV6} concerning a form of a discrete distribution supported at the origin.

{\bf Theorem 2.} {\it
Let $\ae-s=n+\delta, n\in\mathbb N, |\delta|<1/2$. Then a general solution of the equation  \eqref{4} has the following form
\[
\tilde u_d(\xi)=A^{-1}_{d,\neq}(\xi)Q_n(\xi)B_h(Q_n^{-1}(\xi)A^{-1}_{d,=}(\xi)\widetilde{(\ell v_d)}(\xi)) +
 \]
 \[
 +A^{-1}_{d,\neq}(\xi)\left(\sum\limits_{k=0}^{n-1}\tilde c_k(\xi_1)\hat\zeta_2^k+\tilde d_k(\xi_2)\hat\zeta_1^k\right),
\]
where $Q_n(\xi)$ is an arbitrary polynomial of order  $n$ of variables  $\zeta_k=\hbar(e^{-ih\xi_k}-1), k=1,2,$ satisfying the condition
\eqref{100} for $\alpha=n$, $\tilde c_k(\xi_1), \tilde d_k(\xi_2), k=0,1,\cdots,n-1,$ -- are arbitrary functions from  $H^{s_k}(h\mathbb T), s_k=s-\ae+k-1/2$.

The a priori estimate
\[
||u_d||_s\leq const\left(||f||^+_{s-\alpha}+\sum\limits_{k=0}^{n-1}([c_k]_{s_k}+[d_k]_{s_k})\right),
\]
holds, where $[\cdot]_{s_k}$ denotes a norm in  $H^{s_k}(h\mathbb T)$, and
$const$ does not depend on $h$.
}

{\bf Proof.}
We start from the equality  \eqref{5}. Let $Q_n(\xi)$ be an arbitrary polynomial of order $n$ of variables $\zeta_k=\hbar(e^{-ih\xi_k}-1), k=1,2,$ satisfying the condition
\eqref{100} for $\alpha=n$. We multiply the equality  \eqref{5} by  $Q^{-1}_n(\xi)$
\begin{equation}\label{6}
Q^{-1}_n(\xi)A_{d,\neq}(\xi)\tilde u_d(\xi)+Q^{-1}_n(\xi)A^{-1}_{d,=}(\xi)\tilde w_d(\xi)=Q^{-1}_n(\xi)A^{-1}_{d,=}(\xi)\widetilde{(\ell v_d)}(\xi),~~~\xi\in\hbar\mathbb T^2,
\end{equation}

We have in view of Lemma 1
\[
Q^{-1}_n(\xi)A^{-1}_{d,=}(\xi)\widetilde{(\ell v_d)}(\xi)\in\widetilde H^{s-\ae+n}(h\mathbb Z^2),
\]
and since $s-\ae+n=-\delta$ then according to Lemma 3 we write the unique decomposition
\[
Q^{-1}_n(\xi)A^{-1}_{d,=}(\xi)\widetilde{(\ell v_d)}(\xi)=F^+_d(\xi)+F^-_d(\xi),
\]
where
\[
F^+_d(\xi)=B_h(Q^{-1}_n(\xi)A^{-1}_{d,=}(\xi)\widetilde{(\ell v_d)}(\xi)),~~~F^-_d(\xi)=(I-B_h)(Q^{-1}_n(\xi)A^{-1}_{d,=}(\xi)\widetilde{(\ell v_d)}(\xi)).
\]

Taking into account this fact we rewrite the equality  \eqref{6} in the form
\[
A_{d,\neq}(\xi)\tilde u_d(\xi)+A^{-1}_{d,=}(\xi)\tilde w_d(\xi)=Q_n(\xi)F^+_d(\xi)+Q_n(\xi)F^-_d(\xi),
\]
and further,
\[
A_{d,\neq}(\xi)\tilde u_d(\xi)-Q_n(\xi)F^+_d(\xi)=Q_n(\xi)F^-_d(\xi)-A^{-1}_{d,=}(\xi)\tilde w_d(\xi),
\]

Since $F^+_d(\xi)\in\widetilde H^{s-\ae+n}(K_d), F^-_d(\xi)\in\widetilde H^{s-\ae+n}(h\mathbb Z^2\setminus K_d)$ then according to Lemma 1 we conclude  $Q_n(\xi)F^+_d(\xi)\in\widetilde H^{s-\ae}(K_d), Q_n(\xi)F^-_d(\xi)\in\widetilde H^{s-\ae}(h\mathbb Z^2\setminus K_d)$. Applying the inverse discrete Fourier transform we obtain an equality for two discrete distributions. The left hand side vanishes at least under one condition  $\tilde x_1<0$ or $\tilde x_2<0$, and the right hand side vanishes under the condition  .$\tilde x_1>0,\tilde x_2>0$, Thus, it should be a discrete distribution supported on sides of the discrete quadrant $\{(\tilde x_1,\tilde x_2)\in h\mathbb Z^2: \{\tilde x_1>0,\tilde x_2=0\}\cup\{\tilde x_1=0,\tilde x_2>0\}\}$. Using corresponding result from \cite{VV6} we obtain the following form for this discrete distribution
\[
\sum\limits_{k=0}^{n-1}\left(c_k(\tilde x_1)(\Delta_2^{(k)}\delta_d)(\tilde x_2)+d_k(\tilde x_2)(\Delta_1^{(k)}\delta_d)(\tilde x_1)\right),
\]
where all summands should be elements of the space $H^{s-\ae}(h\mathbb Z^2)$.

The left question is how much summands we need in the right-hand side. Counting principle is a very simple because every summand should belong to the space $\widetilde H^{s}(\hbar\mathbb T^2)$.

Let us consider the summand $c_k(\xi_1)\zeta_2^k$. Taking into account that order of $A^{-1}_{d,+}(\xi)$ is $-\ae$ we need to verify the finiteness of the $H^{s-\ae}$-norm for $c_k(\xi_1)\zeta_2^k$. We have
\[
||c_k(\Delta^{(k)}_2\delta_d)||^2_{s-\ae}=\int\limits_{\hbar\mathbb T^2}(1+|\zeta^2|)^{s-\ae}||c_k(\xi_1)\zeta_2^k|^2d\xi
\]
\[
=\int\limits_{\hbar\mathbb T^2}(1+|\zeta^2|)^{s-\ae}||c_k(\xi_1)|^2|\zeta_2^k|^2d\xi\leq a_1\hbar^{2(s-\ae+k+1/2)}\int\limits_{\hbar\mathbb T}|c_k(\xi_1)|^2d\xi_1
\]
\[
\leq a_2\int\limits_{\hbar\mathbb T}(1+|{\zeta_1}^2|)^{s-\ae+k+1/2}|c_k(\xi_1)|^2d\xi_1,
\]
and the constants $a_1, a_2$ do not depend on $h$.

The last summand should be $(n-1)$th because for $n$th summand we obtain a positive growth: for $k=n$ we have $s_n=s-\ae-n+1/2=-n-\delta+n+1/2=-\delta+1/2>0$.

A priori estimates can be obtained in the same way described in \cite{VV6}.
$\blacksquare$

\subsection{The Dirichlet discrete boundary condition}

We consider here first simple case with discrete Dirichlet boundary conditions. We suppose in this section that $\ae-s=1+\delta, |\delta|<1/2, v_d\equiv 0.$ It follows from Theorem 2 that we have the following general solution of the equation \eqref{4}
\begin{equation}\label{7}
 \tilde u_d(\xi)=A^{-1}_{d,\neq}(\xi)(\tilde c_0(\xi_1)+\tilde d_0(\xi_2)),
\end{equation}
where $c_0,d_0\in H^{s-\ae-1/2}(\hbar\mathbb Z)$ are arbitrary functions. To determine uniquely these functions we add the discrete Dirichlet conditions on angle sides
\begin{equation}\label{8}
{u_d}_{|_{\tilde x_1=0}}=f_d(\tilde x_2),~~~~~~~~~~~~~{u_d}_{|_{\tilde x_2=0}}=g_d(\tilde x_1).
\end{equation}

Thus, we have the discrete Dirichlet problem \eqref{4},\eqref{8}.

First, we apply the discrete Fourier transform to discrete conditions \eqref{8} and obtain the following form
\begin{equation}\label{9}
\int\limits_{-\hbar\pi}^{\hbar\pi}\tilde u_d(\xi_1,\xi_2)d\xi_1=\tilde f_d(\xi_2),~~~\int\limits_{-\hbar\pi}^{\hbar\pi}\tilde u_d(\xi_1,\xi_2)d\xi_2=\tilde g_d(\xi_1).
\end{equation}

Substituting \eqref{9} into \eqref{7} we obtain the following relations
\[
\int\limits_{-\hbar\pi}^{\hbar\pi}\tilde u_d(\xi)d\xi_1=\int\limits_{-\hbar\pi}^{\hbar\pi}A^{-1}_{d,\neq}(\xi)\tilde c_0(\xi_1)d\xi_1+\tilde d_0(\xi_2)\int\limits_{-\hbar\pi}^{\hbar\pi}A^{-1}_{d,\neq}(\xi)d\xi_1
\]
\[
\int\limits_{-\hbar\pi}^{\hbar\pi}\tilde u_d(\xi)d\xi_2=\tilde c_0(\xi_1)\int\limits_{-\hbar\pi}^{\hbar\pi}A^{-1}_{d,\neq}(\xi)d\xi_2+\int\limits_{-\hbar\pi}^{\hbar\pi}A^{-1}_{d,\neq}(\xi)\tilde d_0(\xi_2)d\xi_2
\]

Let us denote
\[
\int\limits_{-\hbar\pi}^{\hbar\pi}A^{-1}_{d,\neq}(\xi)d\xi_1\equiv\tilde a_0(\xi_2),~~~\int\limits_{-\hbar\pi}^{\hbar\pi}A^{-1}_{d,\neq}(\xi)d\xi_2\equiv \tilde b_0(\xi_1)
\]
and suppose that $\tilde a_0(\xi_2),\tilde b_0(\xi_1)\neq 0, \forall\xi_1\neq 0,\xi_2\neq 0$.

Therefore, we have the following system of two linear integral equations with respect to two unknown functions $\tilde c_0(\xi_1),\tilde d_0(\xi_2)$
\begin{equation}\label{10}
\left\{
\begin{array}{rcl}
\int\limits_{-\hbar\pi}^{\hbar\pi}M_1(\xi)\tilde c_0(\xi_1)d\xi_1+\tilde d_0(\xi_2)=\tilde F_d(\xi_2)\\
\tilde c_0(\xi_1)+\int\limits_{-\hbar\pi}^{\hbar\pi}M_2(\xi)\tilde d_0(\xi_2)d\xi_2=\tilde G_d(\xi_1),
\end{array}
\right.
\end{equation}
where we have used the following notations
\[
\tilde F_d(\xi_2)=\tilde f_d(\xi_2)\tilde a_0^{-1}(\xi_2),~~~\tilde G_d(\xi_1)=\tilde g_d(\xi_1)\tilde b_0^{-1}(\xi_1),
\]
\[
M_1(\xi)=A^{-1}_{d,\neq}(\xi)\tilde a_0^{-1}(\xi_2),~~~M_2(\xi)=A^{-1}_{d,\neq}(\xi)\tilde b_0^{-1}(\xi_1).
\]

Unique solvability conditions for the system \eqref{10} will be equivalent to unique solvability for the discrete Dirichlet problem \eqref{4},\eqref{8}.

Thus, we obtain the following result.

{\bf Proposition 1.} {\it Let $f_d,g_d\in H^{s-1/2}(\mathbb R_+), s>1/2, v_d\equiv 0$. Then the discrete Dirichlet problem \eqref{4},\eqref{8} is reduced to the equivalent system of linear integral equations \eqref{10}.
}

\subsection{Non-local discrete boundary condition}

We consider here the $\ae-s=1+\delta, |\delta|<1/2$ for the equation \eqref{4} with different boundary conditions, namely
\begin{equation}\label{11}
\begin{array}{rcl}
\sum\limits_{\tilde x_1\in h\mathbb Z_+}u_d(\tilde x_1,\tilde x_2)h=f_d(\tilde x_2),~~~\sum\limits_{\tilde x_2\in h\mathbb Z_+}u_d(\tilde x_1,\tilde x_2)h=g_d(\tilde x_1),\\
\sum\limits_{\tilde x\in h\mathbb Z_{++}}u_d(\tilde x_1,\tilde x_2)h^2=0.
\end{array}
\end{equation}

These additional conditions will help us to determine uniquely the unknown functions $c_0,d_0$ in the solution \eqref{7}.

Indeed, using the discrete Fourier transform we rewrite the conditions \eqref{11} as follows
\begin{equation}\label{12}
\tilde u_d(0,\xi_2)=\tilde f_d(\xi_2),~~~\tilde u_d(\xi_1,0)=\tilde g_d(\xi_1),~~~\tilde u_d(0,0)=0.
\end{equation}

Now we substitute the formulas \eqref{12} into \eqref{7}. The first two equality are
\[
\tilde u_d(0,\xi_2)=A^{-1}_{d,\neq}(0,\xi_2)(\tilde c_0(0)+\tilde d_0(\xi_2))=\tilde f_d(\xi_2),
\]
\[
\tilde u_d(\xi_1,0)=A^{-1}_{d,\neq}(\xi_1,0)(\tilde c_0(\xi_1)+\tilde d_0(0))=\tilde g_d(\xi_1).
\]
It implies the following relations according to the third condition
\[
\tilde f_d(0)=\tilde g_d(0),~~~\text{and from which}~~~\tilde c_0(0)+\tilde d_0(0)=0,~~~\text{and}~~~\tilde c_0(0)=\tilde d_0(0)=0.
\]

Then we have at least formally
\begin{equation}\label{13}
\tilde u_d(\xi)=A^{-1}_{d,\neq}(\xi)\left(A_{d,\neq}(\xi_1,0)\tilde g_d(\xi_1)+A_{d,\neq}(0,\xi_2)\tilde f_d(\xi_2)\right)
\end{equation}

It is left to formulate and to prove exactly the obtained result.

{\bf Theorem 3.} {\it Let $f_d,g_d\in H^{s+1/2}(h\mathbb Z), v_d\equiv 0$. Then the discrete problem \eqref{4},\eqref{11} has unique solution which is given by the formula \eqref{13}.

The a priori estimate
\[
||u_d||_s\leq const(||f_d||_{s+1/2}+||g_d||_{s+1/2})
\]
holds with a const non-depending on $h$
}

{\bf Proof.}
We need to prove the a priori estimate only. Let us consider the first summand
\[
||A^{-1}_{d,\neq}(\xi)A_{d,\neq}(\xi_1,0)\tilde g_d(\xi_1)||^2_s=
\]

\[=\int\limits_{\hbar\mathbb T^2}|A^{-1}_{d,\neq}(\xi_1,\xi_2)A_{d,\neq}(\xi_1,0)\tilde g_d(\xi_1)|^2(1+|\zeta^2|)^sd\xi_1d\xi_2\leq
\]
 \[
\leq C\hbar^{2s}\int\limits_{\hbar\mathbb T^2}|g_d(\xi_1)|^2d\xi\leq C_1\hbar^{2s+1}\int\limits_{-\hbar\pi}^{\hbar\pi}|g_d(\xi_1)|^2d\xi_1\leq
\]
\[
\leq C_2\int\limits_{-\hbar\pi}^{\hbar\pi}|g_d(\xi_1)|^2(1+|\zeta_1^2|)^{s+1/2}d\xi_1=||g_d|^2_{s+1/2}.
\]

The second summand has the same estimate.
$\blacksquare$

\section{A comparison between discrete and continuous solutions}

The continuous analogue of the discrete boundary value problem is the following \cite{V1}.

Let $A$ be a pseudo-differential operator with the symbol $A(\xi), \xi=(\xi_1,\xi_2)$ satisfying the condition
\[
c_1(1+|\xi)^{\alpha}\leq|A(\xi)|\leq c_2(1+|\xi)^{\alpha}.
\]
and admitting the wave factorization with respect to the quadrant $K$ with index $\ae$.

We consider the equation
\begin{equation}\label{14}
(Au)(x)=0,~~~x\in K,
\end{equation}
with the following additional conditions
\begin{equation}\label{15}
\int\limits_{0}^{+\infty}u(x_1,x_2)dx_1=f(x_2),~~~\int\limits_{0}^{+\infty}u(x_1,x_2)dx_2=g(x_1),~~~
\int\limits_{-K}u(x)dx=0.
\end{equation}

A solution of the problem \eqref{14},\eqref{15} is sought in the space $H^s(K)$ \cite{V0} and boundary functions are taken from the space $H^{s+1/2}(\mathbb R_+)$. Such problem was considered in \cite{V1} and it has the solution
\begin{equation}\label{16}
\tilde u(\xi)=A^{-1}_{\neq}(\xi)\left(A_{\neq}(\xi_1,0)\tilde g(\xi_1)+A_{\neq}(0,\xi_2)\tilde f(\xi_2)\right)
\end{equation}
under condition that the symbol $A(\xi)$ admits the wave factorization with respect to the quadrant $K$
\[
A(\xi)=A_{\neq}(\xi)A_=(\xi)
\]
with index $\ae$ such that $\ae-s=1+\delta, |\delta|<1/2$.

To construct a discrete boundary value problem which is good approximation to \eqref{14},\eqref{15} we need to choose $A_d(\xi)$ and $f_d,g_d$ in a special way. First, we introduce the operator $l_h$ which acts as follows. For a function $u$ defined in $\mathbb R$ we take its Fourier transform $\tilde f$ then we take its restriction on $\hbar T$ and periodically extend it to $\mathbb R$. Finally, we take its inverse discrete Fourier transform and obtain the function of discrete variable $(l_hu)(\tilde x), \tilde x\in h\mathbb R$. Thus, we put
\[
f_d=l_hf,~~~g_d=l_hg.
\]

Second, the symbol of digital operator $A_d$ we construct in the same way. If we have the wave factorization for the symbol $A(\xi)$ then we take restrictions of factors on $\hbar\mathbb T^2$ and the periodic symbol $A_d(\xi)$ is a product of these restrictions. For such $f_d,g_d$ and the symbol $A_d(\xi)$ we obtain the following result.

{\bf Theorem 4.}
Let $f,g\in S(\mathbb R), \ae>1$. Then we have the following estimate for solutions $u$ and$u_d$ of the continuous problem \eqref{14},\eqref{15} and the discrete one \eqref{4},\eqref{11}
\[
|u(\tilde x)-u_d(\tilde x)|\leq C(f,g)h^{\beta},
\]
where the const $C(f,g)$ depends on functions $f,g$, $\beta>0$ can be an arbitrary number.

{\bf Proof.}
We need to compare two functions \eqref{13} and \eqref{16}, more exactly their inverse discrete Fourier transform and inverse Fourier transform at points $\tilde x\in K_d$. We have
\[
u_d(\tilde x)-u(\tilde x)=\frac{1}{4\pi^2}\left(\int\limits_{\hbar\mathbb T^2}e^{i\tilde x\cdot\xi}\tilde u_d(\xi)d\xi-\int\limits_{\mathbb R^2}e^{i\tilde x\cdot\xi}\tilde u(\xi)d\xi\right)=
\]
\[
=\frac{1}{4\pi^2}\int\limits_{\mathbb R^2\setminus\hbar\mathbb T^2}e^{i\tilde x\cdot\xi}A^{-1}_{\neq}(\xi)\left(A_{\neq}(\xi_1,0)\tilde g(\xi_1)+A_{\neq}(0,\xi_2)\tilde f(\xi_2)\right)d\xi,
\]
since according to our choice for $A_d,f_d,g_d$ the functions $\tilde u_d$ and $\tilde u$ coincide in points $\xi\in\hbar\mathbb T^2$.

We will estimate one summand.
\[
\left|\frac{1}{4\pi^2}\int\limits_{\mathbb R^2\setminus\hbar\mathbb T^2}e^{i\tilde x\cdot\xi}A^{-1}_{\neq}(\xi)A_{\neq}(\xi_1,0)\tilde g(\xi_1)d\xi\right|\leq
\]
\[
\leq~C\int\limits_{\hbar\pi}^{+\infty}\frac{d\xi_2}{(1+|\xi_1|+|\xi_2|)^{\ae}}\int\limits_{\hbar\pi}^{+\infty}|\xi_1|^{-\gamma}d\xi_1,
\]
since $\tilde g\in S(\mathbb R)$. It implies the required estimate.
$\blacksquare$

\section*{Conclusion}

In this paper we have considered two-dimensional cone only, but the authors continue to work in multidimensional situations and we hope to obtain results similar to a discrete half-space \cite{VV6,TV}.

As first practical applications the authors plan to study discrete variant of a quarter-plane problem in diffraction theory \cite{SFO} and elasticity theory \cite{V0}. We hope it will useful application of the developed theory.

\bibliographystyle{amsplain}
\bibliography{vasilyev}


\end{document}